\documentstyle[12pt]{article}
\begin{document}

\newcommand{\text}[1]{\mbox{{\rm #1}}}  
\newcommand{\gd}{\delta}
\newcommand{\itms}[1]{\item[[#1]]} 
\newcommand{\nin}{\in\!\!\!\!\!/}
\newcommand{\sub}{\subset} 
\newcommand{\cntd}{\subseteq}     
\newcommand{\go}{\omega} 
\newcommand{\Pa}{P_{a^\nu,1}(U)} 
\newcommand{\fx}{f(x)}  
\newcommand{\fy}{f(y)} 
\newcommand{\gD}{\Delta}
\newcommand{\gl}{\lambda} 
\newcommand{\gL}{\Lambda} 
\newcommand{\half}{\frac{1}{2}} 
\newcommand{\sto}[1]{#1^{(1)}}
\newcommand{\stt}[1]{#1^{(2)}}
\newcommand{\Z}{\hbox{\sf Z\kern-0.720em\hbox{ Z}}}
\newcommand{\singcolb}[2]{\left(\begin{array}{c}#1\\#2
\end{array}\right)} 
\newcommand{\ga}{\alpha}
\newcommand{\gb}{\beta} 
\newcommand{\gga}{\gamma}
\newcommand{\ul}{\underline} 
\newcommand{\ol}{\overline} 
\newcommand{\qed}{\kern 5pt\vrule height8pt width6.5pt depth2pt}
\newcommand{\Lrraro}{\Longrightarrow}
\newcommand{\Nb}{|\!\!/}
\newcommand{\NN}{{\rm I\!N}}
\newcommand{\bsl}{\backslash}     
\newcommand{\gt}{\theta}
\newcommand{\op}{\oplus}
\newcommand{\Op}{\bigoplus}          
\newcommand{\CR}{{\cal R}}
\newcommand{\tr}{\bigtriangleup}
\newcommand{\grr}{\omega_1} 
\newcommand{\ben}{\begin{enumerate}}
\newcommand{\een}{\end{enumerate}}
\newcommand{\ndiv}{\not\mid}
\newcommand{\bab}{\bowtie}
\newcommand{\hal}{\leftharpoonup}
\newcommand{\har}{\rightharpoonup}
\newcommand{\ot}{\otimes}
\newcommand{\OT}{\bigotimes}
\newcommand{\bwe}{\bigwedge}
\newcommand{\gep}{\varepsilon}
\newcommand{\gs}{\sigma} 
\newcommand{\rbraces}[1]{\left( #1 \right)}
\newcommand{\bbox}{$\;\;\rule{2mm}{2mm}$}
\newcommand{\sbraces}[1]{\left[ #1 \right]}
\newcommand{\bbraces}[1]{\left\{ #1 \right\}}
\newcommand{\OO}{_{(1)}}
\newcommand{\TT}{_{(2)}}
\newcommand{\FF}{_{(3)}}
\newcommand{\minus}{^{-1}}
\newcommand{\CV}{\cal V} 
\newcommand{\CVs}{\cal{V}_s} 
\newcommand{\un}{U_q(sl_n)'}
\newcommand{\on}{O_q(SL_n)'}
\newcommand{\slq}{U_q(sl_2)}
\newcommand{\olq}{O_q(SL_2)}
\newcommand{\UU}{U_{(N,\nu,\go)}}
\newcommand{\HH}{H_{n,q,N,\nu}} 
\newcommand{\ct}{\centerline}
\newcommand{\bs}{\bigskip}
\newcommand{\qua}{\rm quasitriangular}   
\newcommand{\ms}{\medskip}
\newcommand{\noin}{\noindent}
\newcommand{\mat}[1]{$\;{#1}\;$}
\newcommand{\raro}{\rightarrow}
\newcommand{\map}[3]{{#1}\::\:{#2}\raro{#3}}
\newcommand{\alg}{{\rm Alg}}
\def\newtheorems{\newtheorem{theorem}{Theorem}[subsection]
                 \newtheorem{cor}[theorem]{Corollary}
                 \newtheorem{prop}[theorem]{Proposition}
                 \newtheorem{lemma}[theorem]{Lemma}
                 \newtheorem{defn}[theorem]{Definition}
                 \newtheorem{Theorem}{Theorem}[section]
                 \newtheorem{Corollary}[Theorem]{Corollary}
                 \newtheorem{Proposition}[Theorem]{Proposition}
                 \newtheorem{Lemma}[Theorem]{Lemma}
                 \newtheorem{Defn}[Theorem]{Definition}
                 \newtheorem{Example}[Theorem]{Example}
                 \newtheorem{Remark}[Theorem]{Remark} 
                 \newtheorem{claim}[theorem]{Claim}
                 \newtheorem{sublemma}[theorem]{Sublemma}
                 \newtheorem{example}[theorem]{Example}
                 \newtheorem{remark}[theorem]{Remark}
                 \newtheorem{question}[theorem]{Question}
                 \newtheorem{conjecture}{Conjecture}[subsection]}
\newtheorems
\newcommand{\proof}{\par\noindent{\bf Proof:}\quad}
\newcommand{\dmatr}[2]{\left(\begin{array}{c}{#1}\\
                            {#2}\end{array}\right)}
\newcommand{\doubcolb}[4]{\left(\begin{array}{cc}#1&#2\\
#3&#4\end{array}\right)}
\newcommand{\qmatrl}[4]{\left(\begin{array}{ll}{#1}&{#2}\\
                            {#3}&{#4}\end{array}\right)}
\newcommand{\qmatrc}[4]{\left(\begin{array}{cc}{#1}&{#2}\\
                            {#3}&{#4}\end{array}\right)}
\newcommand{\qmatrr}[4]{\left(\begin{array}{rr}{#1}&{#2}\\
                            {#3}&{#4}\end{array}\right)}
\newcommand{\smatr}[2]{\left(\begin{array}{c}{#1}\\
                            \vdots\\{#2}\end{array}\right)}

\newcommand{\ddet}[2]{\left[\begin{array}{c}{#1}\\
                           {#2}\end{array}\right]}
\newcommand{\qdetl}[4]{\left[\begin{array}{ll}{#1}&{#2}\\
                           {#3}&{#4}\end{array}\right]}
\newcommand{\qdetc}[4]{\left[\begin{array}{cc}{#1}&{#2}\\
                           {#3}&{#4}\end{array}\right]}
\newcommand{\qdetr}[4]{\left[\begin{array}{rr}{#1}&{#2}\\
                           {#3}&{#4}\end{array}\right]}

\newcommand{\qbracl}[4]{\left\{\begin{array}{ll}{#1}&{#2}\\
                           {#3}&{#4}\end{array}\right.}
\newcommand{\qbracr}[4]{\left.\begin{array}{ll}{#1}&{#2}\\
                           {#3}&{#4}\end{array}\right\}}

\title{On Pointed Hopf Algebras and Kaplansky's 10th Conjecture}
\author{Shlomo Gelaki
\\Department of Mathematics\\
Harvard University\\Cambridge, MA 02138}
\date{March 26, 1998}
\maketitle
Radford constructed two families of finite dimensional pointed Hopf
algebras over a field $k$ [R1, 5.1 and 5.2]. The first
one, denoted by \mat{H_{n,q,N,\nu}}, is a family of 
pointed Hopf algebras which contains a sub-family of self dual pointed 
Hopf algebras, denoted by $H_{(N,\nu,\go)}.$ This sub-family generalizes
Taft's well known Hopf algebras [G1,T].
The second one, denoted by \mat{U_{(N,\nu,\go)},}is a family
of pointed, unimodular and ribbon Hopf algebras
constructed for the purpose of computing invariants of knots, links and
3-manifolds.
This family generalizes the well known quantum
group $U_q(sl_2)'$ when $q$ is a root of unity. 
The author constructed a new family of finite dimensional Hopf 
algebras, denoted by ${\cal H}_{n,q,N,\nu,\ga },$ which generalizes 
Radford's \mat{H_{n,q,N,\nu}}, and furthermore proved 
that this new family captures all Hopf
algebras generated as algebras by one grouplike element of order
$N,$ and one non-trivial skew primitive element [G3, 
Theorem 1.1.1]. He also proved that
$U_q(sl_2)'$ and the $8-$dimensional Hopf algebra $U_{(2,1,-1)}$ capture
all minimal quasitriangular Hopf
algebras $U$, generated by one grouplike element of {\em prime} order
$N,$ and two non-trivial $kG(U)$ independent skew primitive elements [G3, 
Theorem 1.2.2]. The purpose of this paper is three fold. First, to 
study the family ${\cal H}_{n,q,N,\nu,\ga },$ and then to construct and 
study a new family of finite dimensional pointed and unimodular Hopf 
algebras which generalizes Radford's \mat{U_{(N,\nu,\go)}.} Second, 
{\bf to answer in the negative Kaplansky's 10th conjecture from 1975} on the
finite number of Hopf algebras of a given dimension \cite{k}.
Third, to
generalize [G3, Theorem 1.2.2] and characterize the 
sub-family \mat{U_{(N,\nu,\go)}} 
via minimal quasitriangularity. 

The paper is organized as follows. In Section 1 we discuss
some aspects of the theory of finite dimensional Hopf algebras, and 
introduce all the Hopf algebras which are relevant for this paper.

In Section 2 we study the family of Hopf algebras ${\cal
H}_{n,q,N,\nu,\ga}$ and its family of dual Hopf algebras ${\cal 
H}_{n,q,N,\nu,\ga}^*.$ 
First, in {\bf Proposition \ref{generalh}} we show that ${\cal
H}_{n,q,N,\nu,\ga}$ is pointed, find a linear basis for it, 
find the first term of its coradical filtration,
find its left 
and right integrals and show that it is an extension of $k\Z_{N/n}$ by 
$H_{n,q,n,\nu}.$  Second, in {\bf Lemma \ref{dualgl}}
we describe $G({\cal H}_{n,q,N,\nu,\ga}^*)$ (where $G(H)$ denotes the
group of grouplike elements of a Hopf algebra $H$), and then
use it in {\bf Theorem \ref{selfdual}} to show that Radford's family 
$H_{(N,\nu,\go)}$ is characterized via 
self duality. Third, in {\bf Theorem \ref{dualpoint}} we find a
necessary and sufficient condition for the dual 
Hopf algebra ${\cal H}_{n,q,N,\nu,\ga}^*$ to be pointed by constructing 
irreducible representations of ${\cal H}_{n,q,N,\nu,\ga}.$ Forth, in {\bf 
Proposition \ref{rad}} we prove that Radford's family $H_{n,q,N,\nu}$ is 
closed under duality, and then use it in {\bf Corollary \ref{new}} to
show that ${\cal H}_{n,q,N,\nu,\ga}^*$ is an extension of
$H_{n,q,n,\nu}^*\cong H_{\frac{n}{(n,\nu)},\go ^\nu,N,\mu}$ 
(hence of $H_{n,q,n,\nu}$ if $\ga\ne 0$) by $k\Z_{N/n}.$ 
Finally, in {\bf Lemma \ref{subhopf}} we find all the Hopf subalgebras of 
$H_{n,q,N,\nu},$ and then use it and Theorem \ref{dualpoint} to determine 
in {\bf Theorem \ref{qt}} when ${\cal H}_{n,q,N,\nu,\ga}$ is quasitriangular.

In Section 3 we first construct and study a new family of finite
dimensional pointed Hopf algebras, denoted by 
${\cal U}_{(n,N,\nu,q,\ga ,\gb ,\gga )}.$ In
particular in {\bf Proposition \ref{generalu}} we show that it is 
pointed, find a linear basis for it, find the first term of its coradical 
filtration, show that it is an extension of
$k\Z_{N/n}$ by $U_{(n,\nu,\go )},$ for some primitive $nth$ root of unity 
$\go ,$ and prove that it is always unimodular by 
computing a non-zero $2-$sided integral for it. We also show that over any
infinite field which contains a primitive $nth$ root of unity, 
our new family contains infinitely many non-isomorphic Hopf
algebras of any dimension of the form $Nn^2,$ where $2<n<N$ are 
integers so that $n$ divides $N.$ Thus, in {\bf Corollary \ref{kap}} 
we prove that Kaplansky's 10th conjecture from 1975 [K] is false. Second, in {\bf 
Theorem \ref{u}} we use Theorem \ref{dualpoint} to generalize [G3, Theorem
1.2.2], and show that over an algebraically
closed field of characteristic zero, Radford's family $U_{(N,\nu,\go)},$ 
which we identify as the sub-family ${\cal  U}_{(\frac{N}{(N,\nu)},N,\nu,
\go ^{\nu},0,0,1)},$  captures all minimal quasitriangular Hopf 
algebras $U$, generated as algebras by one grouplike element of an
arbitrary order
$N,$ and two non-trivial $kG(U)$ independent skew primitive elements.  

I wish to thank the referee for his/her very useful comments.
\section{Preliminaries}
Throughout this paper, unless otherwise explicitly stated, $k$ denotes 
a field. Also, we shall denote the greatest common divisor of two 
integers $a$ and $b$ by $(a,b).$ 

The Hopf algebras which are studied in this paper are pointed. 
Recall that a Hopf algebra $H$ over the field $k$ is pointed if its simple
subcoalgebras are
$1-$dimensional; that is, they are generated by elements of $G(H).$ For any 
$g,h\in G(H)$ we denote the vector space of $g:h$ skew primitives of 
$H$ by $P_{g,h}(H)=\{x\in H\:|\gD(x)=x\ot g+h\ot x\}.$
Thus the classical primitive elements are $P_{1,1}(H)$.
The element $g-h$ is always $g:h$ skew primitive. If $P_{g,h}(H)$ consists 
only of the linear span of $g-h$ we say it is trivial. Otherwise, let  
$P_{g,h}(H)^{'}$ denote a complement of $sp_k\{g-h\}$ in  
$P_{g,h}(H)$. Taft-Wilson theorem [TW] states that  
the first term $H_1$ of the coradical filtration of $H$ is given by:
\begin{equation}\label{wilson} 
H_1=kG(H)\bigoplus (\Op_{g,h\in G(H)} P_{g,h}(H)^{'}).
\end{equation}

The following lemma is very useful in determining whether certain Hopf 
algebras are pointed. Various forms of it appear in the literature; see
[R1, Lemma 1] for example. 
\begin{Lemma}\label{point}
Suppose that $H$ is a Hopf algebra over the field $k$ which is generated
as an algebra by a subset $S$ of $G(H)$ and by $g:g'$ skew primitives, where 
$g,g'$ run over $S.$ Then $H$ is pointed and $G(H)$ is generated as a 
group by $S.$ 
\end{Lemma}

The next proposition is very useful, in particular in the context of
quantum groups; see [R1, Proposition 1] for example.
\begin{Proposition}\label{qbinom}
Suppose that $A$ is an algebra over the field $k$ and $\go\in k$ is a
primitive $nth$ root of unity. Let \mat{a,x\in A} satisfy \mat{xa=\go a x.} 
Then $(a+x)^n=a^n+x^n.$
\end{Proposition} 
\begin{Lemma} {\bf [G2, Lemma 0.2]}\label{trivial}
Let $H$ be a finite dimensional Hopf algebra 
over a 
field $k.$ Assume that $(dimH)1\ne 0$ and let $x\in P_{a,1}(H)$. Then 
$xa=ax$ if and only if $x\in sp_k \{a-1\}$.
\end{Lemma}

Next, we give a very brief description
of the underlying Hopf algebra structure of $D(H)$, the
Drinfel'd double of a finite dimensional Hopf algebra $H$ [D]. As a
coalgebra $D(H)=H^{*cop}\ot H$ (where
$H^{*cop}=H^*$ as an algebra, with new comultiplication $\gD^{cop}=\tau
\circ \gD ,$ where $\tau:H^*\ot H^*\raro H^*\ot H^*$ is the usual twist
map). Set \mat{p\bab a =p\ot a} for $p\in H^*$ and $a\in H$.
Multiplication
in the double can be described in several ways. The following are very
convenient:
$$(p\bab a)(q\bab b)=\sum p\rbraces{a\OO\cdot q\cdot s\minus(a\FF)}\bab a\TT b$$
and
$$(p\bab a)(q\bab b)=\sum pq\TT\bab\rbraces{(s^*)\minus(q\OO)
\har a \hal q\FF}b$$
where $<a\cdot q,b>=<q,ba>,$ $<q\cdot a,b>=<q,ab>,$ $p\har a=\sum 
a\OO<p,a\TT>$ and $a\hal p=\sum<p,a\OO>a\TT$ for all $p,q\in H^*$
and $a,b\in H.$ 
The antipode $S$ of $D(H)$ is given
by {\rm S}$(p\bab a)=(\gep\bab s(a))((s^*)\minus(p)\bab1)$ for $p\in H^*$
and $a\in H$. Observe that \mat{(p\bab 1)(q\bab b)=pq\bab b} and 
\mat{(p\bab a)(\gep\bab b)=p\bab ab.} As a consequence, the one-one maps
$\map{i_H}{H}{D(H)}$ and $\map{i_{H^*}}{H^{*cop}}{D(H)}$ defined by
\mat{i_H(h)=\gep\bab h} and \mat{i_{H^*}(p)=p\bab1} respectively for
$h\in H$ and $p\in H^*$ are Hopf algebra maps.

We recall now the definition of a quasitriangular Hopf
algebra. Let $H$ be a Hopf algebra and $R=\sum R^{(1)}\ot R^{(2)}\in H\ot
H.$ 
The pair $(H,R)$ is said to be a quasitriangular Hopf algebra
if the following axioms hold ($r=R$):
\begin{description}
\item[(QT.1)] \mat{\sum\gD(\sto{R})\ot\stt{R}=
              \sum\sto{R}\ot\sto{r}\ot\stt{R}\stt{r} },
\item[(QT.2)] \mat{\sum\gep(\sto{R})\stt{R}=1 },
\item[(QT.3)] \mat{\sum\sto{R}\ot\gD^{cop}(\stt{R})=
              \sum\sto{R}\sto{r}\ot\stt{R}\ot\stt{r} },
\item[(QT.4)] \mat{\sum\sto{R}\gep(\stt{R})=1} and
\item[(QT.5)] \mat{\rbraces{\gD^{cop}(h)}R=R(\gD(h))} for all $h\in H.$
\end{description}
Drinfel'd has shown that the double of any finite dimensional Hopf algebra
is quasitriangular [D]. Suppose that \mat{\map{f}{H}{H'}} is a
surjective Hopf algebra map and set \mat{R'=\rbraces{f\ot f}(R).} Then
$(H',R')$ is quasitriangular. Let $(H,R)$ be quasitriangular and write 
$R=\sum\sto{R}\ot\stt{R}$ in the shortest possible way. Set
$l(R)=sp_k\bbraces{\sto{R}}$ and $r(R)=sp_k
\bbraces{\stt{R}}.$ Let $H_R$ be the Hopf subalgebra of $H$ generated
by $l(R)$ and $r(R)$. Then $l(R)$ and $r(R)$
are finite dimensional Hopf subalgebras of $H,$ and there exists an
isomorphism of
Hopf algebras \mat{f:l(R)^{*cop}{\raro}r(R)} and  a
unique  surjection of Hopf algebras $F:D(l(R))\raro H_{R}$ which 
satisfies \mat{
F_{\mid_{l(R)}}=i(l(R))} and \mat{F_{\mid_{l(R)}*cop}=f} [R2].
If $H=H_R$ then $(H,R)$ is called a {\em minimal} 
quasitriangular Hopf algebra. We shall also say that $H$ is a (minimal) 
quasitriangular Hopf algebra if there exists $R\in H\ot H$ such that 
($H,R$) is a (minimal) quasitriangular Hopf algebra.       

We now recall the definition of $H_{n,q,N,\nu}$ and $H_{(N,\nu,\go)}$
[R1, 5.1].
Let $n,\nu$ and $N$ be positive integers such that $n|N$
and $1\leq\nu<N$. Suppose $q\in k$ is a primitive $nth$ root of unity and 
let $r=|q^{\nu}|$. As an algebra $H=H_{n,q,N,\nu}$ is generated by a and
$x$ which satisfy the relations
$$a^N=1,\ \ \ \ x^r=0\;\;{\rm and}\;\;\;\;\; xa=qax.$$
The coalgebra structure of $H$ is determined by 
$$\Delta(a)=a\otimes a,\;\;\Delta (x)=x\otimes a^\nu+1\otimes x,\;\;
\gep(a)=1 \;\;{\rm and}\;\; \gep(x)=0.$$     
The antipode of $H$ is determined by
$$s(a)=a\minus\;\;\;\;{\rm and}\;\;\;\; s(x)=-q^{-\nu}a^{-\nu}x.$$
By [R1, Proposition 8 and Theorem 4 a)] and Lemma
\ref{point}, $H$ satisfies the following: 
\begin{description}
\item [(H.1)] As an algebra $H$ is generated by a grouplike element $a$ 
      of order $N,$ and by the space $P_{a^\nu,1}(H)$ of dimension$\le 2.$
\item [(H.2)] $H$ is pointed with $G(H)=(a),$ and if $r>1$ then it is 
      self dual if and only if there exists a
      primitive $Nth$ root of unity $\go \in k$ such that $q=\go ^\nu.$
\end{description}
If the field $k$ contains a primitive $Nth$ root of unity $\go $
such that $q=\go ^\nu$ and $N\Nb \nu ^2,$ then $H_{n,q,N,\nu}$ is
denoted by $H_{(N,\nu,\go )}.$ 

We next recall the definition of the family ${\cal 
H}_{n,q,N,\nu,\ga}$ [G3, 1.1].
Let $n,\nu$ and $N$ be positive integers such that $n|N$
and $1\leq\nu<N,$ and let $\ga\in k.$ Suppose $q\in k$ is a primitive 
$nth$ root of unity and let $r=|q^{\nu}|$. For $\ga =0$ we define
${\cal H}_0={\cal H}_{n,q,N,\nu,0}= H_{n,q,N,\nu}.$ 
For {\bf $0\ne\ga\in k$} we further assume 
that {\bf $(n,\nu)=1$} (hence $r=n$), $n>1$ and {\bf $N\Nb\nu n.$} In this
case, 
as an algebra ${\cal H}_{\ga }={\cal H}_{n,q,N,\nu,\ga}$ is generated by a
and $x$ which satisfy the relations
$$a^N=1,\;x^n=\ga(a^{\nu n}-1)\;\;and\;\;xa=qax.$$
The coalgebra structure of ${\cal H}_{\ga }$ is determined by 
$$\Delta(a)=a\otimes a,\;\;\Delta (x)=x\otimes a^\nu+1\otimes x,\;\;
\gep(a)=1 \;\;{\rm and}\;\; \gep(x)=0.$$ 
The antipode of ${\cal H}_{\ga }$ is determined by 
$$s(a)=a\minus\;\;\;\;{\rm and}\;\;\;\; s(x)=-q^{-\nu}a^{-\nu}x.$$
Note that since $N\Nb\nu n$ it follows that $x^n\ne 0.$
\begin{Remark}\label{radf1}
{\rm Suppose $\ga \ne 0.$ Then to construct ${\cal H}_{\ga }={\cal 
H}_{n,q,N,\nu,\ga}$ we follow 
the line of argument and notations used by Radford in order to construct 
$H_{n,q,N,\nu}$ [R1, 5.1]. We let $C$ be the $4-$dimensional 
coalgebra 
over $k$ with basis $\{A,B,D,X\}$ whose structure is determined by
\begin{equation}\label{structure}
\gD(A)=A\ot A,\;\;\gD(B)=B\ot B,\;\;\gD(D)=D\ot D\;\;and\;\;\gD(X)=X\ot B+ 
D\ot X,
\end{equation}
and let $I$ be the ideal of the tensor algebra $T(C)$ generated by
$$A^N-1,\;\;B-A^{\nu},\;\;D-1,\;\;XA-qAX\;\;and\;\;X^n-\ga (A^{\nu n}-1).$$
It is straightforward to check that $T(C)$ has a unique bialgebra structure 
determined by (\ref{structure}), that $I$ is a biideal of $T(C)$ and 
finally that ${\cal H}_{\ga }=T(C)/I,$ where $a=\overline{A}$ and 
$x=\overline{X},$ is a Hopf algebra.} 
\end{Remark}

In the following we state that over an algebraically
closed field of characteristic zero, the family
${\cal H}_{n,q,N,\nu,\ga}$ is unique with respect to property (H.1).
We also prove it here for the sake of completeness. 
\begin{Theorem} {\bf [G3, Theorem 1.1.1]}\label{h} 
Let $k$ be an algebraically closed field of characteristic zero. Let 
${\cal H}$ be a Hopf algebra over $k,$ and suppose
there exist integers $1\leq\nu<N$ so that 
as an algebra ${\cal H}$ is generated by a grouplike element $a$ of order 
$N,$ and by the space $P_{a^\nu,1}({\cal H})$ of dimension$\le 2.$
Then there exist a positive integer $n$ such that $n|N,$ $\ga \in k$ and
a primitive $nth$ root of unity $q\in k$ such that 
${\cal H}={\cal H}_{n,q,N,\nu, \ga}$.
\end{Theorem}
\proof 
If $dim_kP_{a^\nu,1}({\cal H})=1$ then $P_{a^\nu,1}({\cal H})=sp_k\{a^\nu 
-1\}.$ Thus, ${\cal H}=k(a)={\cal H}_{1,1,N,\nu,0},$ and we are done. Now 
assume $dim_kP_{a^\nu,1}({\cal H})=2.$  First note that the semisimple 
commutative
algebra $k(a)$ acts on $P_{a^\nu,1}({\cal H})$ by conjugation and 
that $sp_k\bbraces{a^\nu-1}$ is a $k(a)-$submodule
of $P_{a^\nu,1}({\cal H})$. Thus there exists a $1-$dimensional
$k(a)$ complement of it. Let $\{x\}$ be a $k$ basis of this 
complement. Then:
\begin{equation}\label{two}
a\minus xa=qx\;\;\mbox{or equivalently}\;\;xa=qax\;\;\mbox{for 
some}\;\;q\in k.
\end{equation}
Since $a^N=1$ it follows that $q^N=1$. Set $n=|q|$. Then $n|N$. Also set 
$r=|q^\nu|$.
By (\ref{two}) $(x\ot a^\nu)(1\ot x)=q^{-\nu}(1\ot x)(x\ot a^\nu)$,
hence by Proposition \ref{qbinom}  
$$\Delta (x^r)=(x\ot a^\nu+1\ot x)^r=x^r\ot a^{\nu r}+1\ot x^r.$$
Moreover,
$$x^r a^{\nu r}=q^{\nu r^2} a^{\nu r}x^r= a^{\nu r}x^r.$$
Thus it follows from Lemma \ref{trivial} that $x^r=\ga(a^{\nu r}-1)$ for 
some $\ga\in k$. Therefore,
if either $\ga=0$ or $N|\nu r$ then there exists a surjection of Hopf 
algebras ${\cal H}_{n,q,N,\nu, 0}\raro {\cal H},$ which is also an 
injection by Proposition 2.1(3) and [M, Theorem 5.3.1], and we are done.
If $\ga \ne 0$ and $N\Nb \nu r$ then $x^r\ne 0.$ Since
$\ga(a^{\nu r}-1)$ and $a$
commute, so do $x^r$ and $a.$ But, $x^r\ne 0,$ hence $q^r=1$
which is equivalent to 
$n|r=\frac{n}{(n,\nu)}.$ Thus, $(n,\nu)=1$ and $r=n,$ and 
we conclude that there exists a surjection of Hopf
algebras ${\cal H}_{n,q,N,\nu, \ga}\raro {\cal H},$ which is also an
injection by Proposition 2.1(3) and [M, Theorem 5.3.1].\qed

We conclude this section by recalling the definition of $U_{(N,\nu,\go)}$ 
[R1, 
5.2]. Let $\nu$ and $N$ be positive integers such that $1\leq\nu <N$ and
$N\Nb\nu^2$. Suppose that $\go\in k$ is a primitive $Nth$ root
of unity. Let $q=\go^\nu$ and $r=|q^\nu|=|\go^{\nu^2}|$.
As an algebra $U=U_{(N,\nu,\go)}$ is generated by $a,x$ and $y$ which
satisfy the relations
$$a^N=1,\;\;x^r=0,\;\;y^r=0,\;\;xa=qax,\;\;ya=q\minus ay\;\;{\rm and}\;\;
yx-q^{-\nu}xy=a^{2\nu}-1.$$
The coalgebra structure of $U$ is determined by
$$\gD(a)=a\ot a,\;\;\gD(x)=x\ot a^\nu+1\ot x,\;\;\gD(y)=y\ot a^\nu+1\ot y,$$
$$\gep(a)=1,\;\;\;\gep(x)=0 \;\;\;{\rm and}\;\;\;\gep(y)=0.$$
The antipode of $U$ is determined by
$$s(a)=a\minus,\;\;\;s(x)=-q^{-\nu}a^{-\nu}x\;\;\;{\rm and}\;\;\;
s(y)=-q^\nu a^{-\nu} y.$$
It follows from Lemma \ref{point}, [G2, Lemma 1.2.1]
and [R1, Proposition 10 d) and Theorem 5] that $U$
satisfies the following:
\begin{description}
\item [(U.1)] As an algebra $U$ is generated by a grouplike 
      element $a$ of order $N,$
      and by the space $P_{a^\nu,1}(U)=sp\{a^v-1,x,y\}$ of dimension $3.$
      In particular, $G(U)=(a).$
\item [(U.2)] $P_{a^v,1}(U)=sp_k\bbraces{a^v-1}$ for all $v\neq\nu.$
              In particular, it follows from (\ref{wilson}) that  
              $U_1=sp\{a^i,a^jx,a^ly|0\le 
              i,j,l\le N-1\}.$
\item [(U.3)] $U$ is minimal quasitriangular.
\item [(U.4)] $U$ is pointed, ribbon and unimodular.
\end{description}
The author proved that over an algebraically closed field of
characteristic zero, if $U$ is a Hopf algebra which satisfies (U.1)-(U.3) 
with prime $N$ then either $U=U_{(2,1,-1)}$ or
there exists a primitive $Nth$ root 
of unity $q$ so that $U=U_q(sl_2)^\prime$ 
[G3, Theorem 1.2.2]. 
\section{The Family ${\cal H}_{n,q,N,\nu,\ga}$}
In this section we study some properties of the Hopf algebra ${\cal 
H}_{n,q,N,\nu,\ga}$ and of its dual Hopf algebra 
${\cal H}_{n,q,N,\nu,\ga}^*.$ Beside of being interesting on their own 
right they will prove to be crucial in Section 3.
\begin{Proposition}\label{generalh}
Let ${\cal H}_{\ga }={\cal H}_{n,q,N,\nu,\ga}.$ Then:
\ben
\item ${\cal H}_{\ga }$ is a pointed Hopf algebra with $G({\cal H}_{\ga 
      })=(a)$ a cyclic group of order $N.$
\item The set $\{a^ix^j|0\le i \le N-1,0\le j \le r-1\}$ forms a linear 
      basis of ${\cal H}_{\ga },$ hence $dim{\cal H}_{\ga }=Nr.$
\item If $n|\nu$ then 
      $P_{a^\nu,1}({\cal H}_0)=sp_k\{a^\nu-1\}.$
      If $n\Nb\nu$ then $P_{a^\nu,1}({\cal H}_{\ga })=sp_k\{a^\nu-1,x\},$
      $P_{a^i,1}({\cal H}_{\ga })=sp_k\{a^i-1\}$ for all $i\ne \nu,$ 
      and $({\cal H}_{\ga })_1=sp\{a^i,a^jx|0\le i,j\le N-1\}.$
\item The elements $\lambda_l=\left (\frac{1}{N}\sum_{i=0}^{N-1}a^i\right 
      )x^{r-1}$ and $\lambda_r=x^{r-1}\left (\frac{1}{N}\sum_{i=0}^{N-1}a^i
      \right)$ are non-zero left and right integral of ${\cal H}_{\ga}.$ 
      In particular ${\cal H}_{\ga}$ is not unimodular when $r>1,$ and the 
      distinguished grouplike element of ${\cal H}_{\ga}^*,$ $\psi,$ is
      determined by $\psi(a)=q^{r-1}$ and $\psi(x)=0.$
\item As an algebra ${\cal H}_{\ga}\cong
      k\Z_{N/n}\#_{\sigma}H_{n,q,n,\nu}$
      is a crossed product of a group algebra and a Radford Hopf algebra 
      with the trivial action.
\item If the field $k$ is algebraically closed then ${\cal H}_{\ga}\cong 
      {\cal H}_{\gb}$ as Hopf algebras for any $\ga ,\gb \in k^*.$
\een
\end{Proposition}
\proof
1. Follows from Lemma \ref{point}.\\
2. For $\ga =0$ it was proved in [R1, Proposition 7 b)] using 
Bergman's Diamond Lemma (see e.g. [R1, page 214] for a description 
of this lemma). The proof for $\ga \ne 0$ is very similar, and we merely 
sketch it following the line of argument of the proof of 
[R1, Proposition 7 b)]. Using the notation of Remark \ref{radf1},
let $T(C)=k\{A,B,D,X\}$ where $A<B<D<X$ and the substitution rules are:
$$A^N\leftarrow 1,\;\;B\leftarrow A^{\nu},\;\;D\leftarrow 1,\;\;XA\leftarrow 
qAX\;\;and\;\;X^n \leftarrow \ga (A^{\nu n}-1).$$
The only ambiguities to consider are overlap ambiguities, and they are all 
easily resolved. We just note that in the case $\ga \ne 0,$ it is 
necessary to 
have $r=n$ (so that $q^r=q^n=1$) when resolving, for instance, the overlap 
ambiguity $(X^{n-1}X)A=X^{n-1}(XA).$\\ 
3. For $\ga =0$ it was proved in [R1, Theorem 4 a)]. Suppose 
$\ga \ne 0$ (so in particular $n\Nb\nu$), and let $z\in P_{a^i,1}({\cal 
H}_{\ga }).$ Then the result follows in a straightforward manner after 
expressing $z$ in terms of the basis given in part 2, 
to compare the expressions of $\gD(z)$ and $z\ot a^i + 1\ot z.$ 
The description of $({\cal H}_{\ga })_1$ follows now from (\ref{wilson}).\\
4. We show that $\lambda_l$ is a left integral. For $\ga =0$ it was shown 
in [R1, Proposition 7 d)]. Suppose $\ga \ne 0.$ Since $\frac{1}{N}\left 
(\sum_{i=0}^{N-1}a^i\right )$ is a $2-$sided integral of $k(a)$ it follows that 
$a\lambda_l=\lambda_l=\varepsilon(a)\lambda_l.$ 
Now, let $\phi \in G({\cal H}_{\ga}^*)$ be determined by $\phi(a)=q$ and
$\phi(x)=0.$ Using the facts that $q^n=1,$ and $\har$ is an algebra map we 
compute, 
\begin{eqnarray*}
\lefteqn{ 
x\lambda_l=x\left (\frac{1}{N}\sum_{i=0}^{N-1}a^i\right)x^{n-1}}\\
& & =\left (\frac{1}{N}\sum_{i=0}^{N-1}q^ia^i\right)x^{n}=
\left (\frac{1}{N}\sum_{i=0}^{N-1}q^ia^i\right)\ga (a^{\nu n}-1)\\
& & =\left(\phi\har \frac{1}{N}\sum_{i=0}^{N-1}a^i\right)\ga (a^{\nu
n}-1)=
\left (\phi\har \frac{1}{N}\sum_{i=0}^{N-1}a^i\right)(\phi\har (\ga
(a^{\nu n}-1)))\\
& & =\phi\har \left (\left (\frac{1}{N}\sum_{i=0}^{N-1}a^i\right)(\ga
(a^{\nu n}-1))\right)=
\phi\har 0=0=\varepsilon(x)\lambda_l.
\end{eqnarray*}
Similarly, $\lambda_r$ is a right integral of ${\cal H}_\ga.$
Since $\lambda_la=q^{n-1}\lambda_l=q^{-1}\lambda_l$ and $\lambda_lx=0,$ 
it follows that 
$\psi$ is the distinguished grouplike element of ${\cal H}_\ga ^*$ and that 
${\cal H}_\ga $ is not unimodular. \\
5. Since $a^n$ is central in ${\cal H}_\ga $ it follows that
$K=k(a^n)\cong k\Z_{N/n}$ 
is a central Hopf subalgebra of ${\cal H}_\ga ,$ hence is normal. 
Let $\pi :{\cal H}_\ga \raro H_{n,q,n,\nu}$ be the Hopf algebra 
surjection determined by $\pi (a)=a$ and $\pi(x)=x.$ Then ${\cal
H}_\ga K^+=ker\pi,$ hence
the following is a short exact sequence of Hopf algebra maps:
$$K\stackrel{i}{\hookrightarrow}{\cal H}_\ga 
\stackrel{\pi}{\raro}H_{n,q,n,\nu}$$
where $i$ is the inclusion map. Thus, ${\cal H}_{\ga}\cong 
k\Z_{N/n}\#_{\sigma}H_{n,q,n,\nu}.$ Finally, by [DT], the action is the 
trivial one since $k\Z_{N/n}$ is central.\\
6. Let ${\cal H}_\ga =k<a,x>$ and ${\cal H}_\gb =k<g,y>.$ Then it is
straightforward to check that the map $f:{\cal H}_\ga \raro 
{\cal H}_\gb $ determined by $f(a)=g$ and $f(x)=(\ga /\gb )^{1/n}y$ is
an isomorphism of Hopf algebras. This completes the proof of the 
proposition.\qed

We now consider the family of the dual Hopf algebras ${\cal
H}_{n,q,N,\nu,\ga }^*.$ 
\begin{Lemma}\label{dualgl}
Let ${\cal H}_{\ga}={\cal H}_{n,q,N,\nu,\ga},$ and assume $k$ 
contains a primitive $Nth$ root of unity $\go .$ If 
$\ga =0$ then $|G({\cal H}_0^*)|=N,$ and if $\ga \ne 0$ then $|G({\cal 
H}_{\ga}^*)|=\#\{0\le i\le N-1|N\;divides\;\nu ni\}.$ 
\end{Lemma}
\proof
For $\ga =0$ and $n=1$ the statement is clear since ${\cal
H}_{\ga}=k\Z_N.$  
Suppose $\ga =0$ and $n>1,$ and let $\phi \in G({\cal H}_{0}^*).$
Since $\phi$ is an algebra homomorphism and $a^N=1,$ it follows that
$\phi(a)=\go ^i$ for some $0\le i\le N-1,$
and the relation $xa=qax$ yields that $\phi(x)=0.$
Now, it is straightforward to check that the map
$\phi :{\cal H}_{0}\raro k$ determined by $\phi(a)=\go ^i$ for some
$0\le i\le N-1$ and $\phi(x)=0$ is indeed an algebra homomorphism.
Thus, the result follows.

Suppose $\ga \ne 0$ (so in particular $n>1$), and let $\phi \in G({\cal 
H}_{\ga}^*).$  
As above, $\phi(a)=\go ^i$ for some $0\le i\le N-1,$ and $\phi(x)=0.$
Now, since $x^n=\ga (a^{\nu n}-1)$ it follows that
$0=\phi(x)^n=\phi(x^n)=\ga 
\phi(a^{\nu n}-1)=\ga (\go ^{\nu ni}-1)$ which is equivalent to $N|\nu ni.$ 
Therefore, $G({\cal H}_{\ga}^*)$ and $\{0\le i\le N-1|N\;divides\;\nu
ni\}$ are in one to one correspondence and the result follows.\qed
\begin{Remark}\label{approx}
{\rm Since $\{0\le i\le N-1|N\;divides\;\nu ni\}$ is a sub group of 
$\Z_N$ it follows that $\#\{0\le 
i\le N-1|N\;divides\;\nu ni\}$ divides $N.$ Moreover, since
$N/n\le N-1$ is an integer we have that $n\le |G({\cal H}_{\ga}^*)|\le N.$}
\end{Remark}

We are ready now to characterize Radford's $H_{(N,\nu,\go )}$ via self
duality. 
\begin{Theorem}\label{selfdual}
Let ${\cal H}_{\ga}={\cal H}_{n,q,N,\nu,\ga}.$ Then 
${\cal H}_{\ga}$ is self dual if and only if $\ga =0$ and 
either ${\cal H}_0=k\Z_N$ or ${\cal H}_0=H_{(N,\nu,\go )}.$
\end{Theorem}
\proof
The 'only if' part is [R1, Proposition 8]. Conversely, suppose
${\cal H}_{\ga}$ is self dual, and assume on the contrary that $\ga \ne 0.$
Then $|G({\cal H}_{\ga}^*)|=|G({\cal H}_{\ga})|=N,$ 
hence by Lemma \ref{dualgl}, $\#\{0\le i\le N-1|N\;divides\;\nu ni\}=N.$ 
In particular $N|\nu n$ which is a contradiction. Therefore, $\ga =0$ and
the result follows from [R1, Proposition 8].\qed 

We now determine when ${\cal H}_{\ga}^*$ is pointed.
\begin{Theorem}\label{dualpoint}
Let ${\cal H}_{\ga}={\cal H}_{n,q,N,\nu,\ga },$ and assume $k$ is
algebraically closed. Then ${\cal H}_{\ga}^*$ is pointed if and only if
$\ga =0.$ Moreover, if $\ga \ne 0$ then ${\cal H}_{\ga}^*$ contains an
$n^2-$dimensional simple sub coalgebra.
\end{Theorem}
\proof
Suppose $\ga =0.$ If $r=1$ then ${\cal H}_0=k\Z_N,$ hence ${\cal 
H}_0^*=k\Z_N$ is pointed. Suppose $r>1$ (so ${\cal H}_0\ne k\Z_N$). We 
show that ${\cal 
H}_{0}^*$ is pointed by proving that all its irreducible comodules
are $1-$dimensional, and hence that all its simple subcoalgebras are
$1-$dimensional which is to say that it is pointed. 
Equivalently, we prove that all the irreducible ${\cal H}_{0}-$modules 
are $1-$dimensional. Indeed, suppose $V$ is an
irreducible ${\cal H}_{0}-$module so that $dimV>1.$ Let $0\ne v\in V$ so
that $a\cdot v=\lambda v,$ $\lambda\in k.$ Since $a\cdot(x^i\cdot v)=
q^ix^i\cdot(a\cdot v)=q^i\lambda (x^i\cdot v)$ it follows that either
$x^i\cdot v=0$ or $x^i\cdot v$ is an eigenvector of $a,$ all $i.$ Let
$i\ge 0$ be
the minimal integer so that $x^i\cdot v\ne 0$ but $x^{i+1}\cdot v=0$
(such an $i$ exists since $x^r=0$). Then, $sp_k\{x^i\cdot v\}$ is stable
under the action of $a$ and $x,$ hence is a non-trivial sub-module of $V$
which is a contradiction, and the result follows. 

Conversely, suppose $\ga
\ne 0$ (so in particular $n>1$). We show that ${\cal H}_{\ga}^*$ is not 
pointed by constructing
irreducible ${\cal H}_{\ga}-$modules of dimension greater than 1.
It will follow that ${\cal H}_{\ga}^*$ contains simple
subcoalgebras of dimension greater than 1, and hence is not pointed. 
Pick an $Nth$ primitive root of unity $\go \in k,$ and define the 
${\cal H}_{\ga}-$module $V_\go $ as follows: $V_\go $ has a linear basis 
$\{v_0,\dots,v_{n-1}\}$ and the action, $\cdot,$ is determined by
$$a\cdot v_i=q^{-i}\go v_i\;\;\;and\;\;\;x\cdot v_i=\gb v_{i+1}$$
where $\gb =[\ga (\go ^{\nu n}-1)]^{1/n}.$ Note that since $\go \ne 1$ and 
$N$ does not divide $\nu n$ it follows that $\gb \ne 0.$ It is 
straightforward to check that this definition is compatible with the 
defining relations of ${\cal H}_{\ga}.$ For example, $x^n\cdot v_i=\gb 
^nv_{i+n}=\ga (\go ^{\nu n}-1)v_{i}$ and $\ga (a^{\nu n}-1)\cdot 
v_i=\ga ((q^{-i}\go 
)^{\nu n}-1)v_{i}=\ga (\go ^{\nu n}-1)v_{i}$ as $q^n=1.$ We claim that 
$V_\go $ is irreducible. Indeed, let $0\ne V\subseteq V_{\go }$ be a 
submodule. Since $a\cdot:V\raro V$ is diagonalizable, $V_\go 
=\oplus_{i=0}^{n-1}sp_k\{v_{i}\}$ and $v_{i},v_{j}$ are eigenvectors of
distinct eigenvalues for $0\le i\ne j\le n-1,$
it follows that there exists $0\le i\le 
n-1$ so that $v_{i}\in V.$ But then, $v_{j}=x^{j-i}\cdot v_i\in V$ for
all $j,$ and hence $V=V_\go .$ Therefore, if ${\cal H}_{\ga}^*$ is pointed
then $\ga=0.$ This concludes the proof of
the theorem.\qed  

By [R1, Proposition 8], $H_{(N,\nu,\go)}$ is self dual, so in particular its 
dual is pointed. By Theorem \ref{dualpoint}, $H_{n,q,N,\nu}^*$ is 
always pointed regardless being self dual. In fact we can say more. 
\begin{Proposition}\label{rad}
Let $H=H_{n,q,N,\nu}\ne k\Z_N,$ and let $\go \in k$ be a primitive $Nth$ 
root of 
unity. Then $H^*=H_{\frac{N}{(N,\nu)},\go ^\nu,N,\mu}$ for some 
$0\le\mu\le N-1$ so that $q=\go ^\mu.$ 
\end{Proposition}
\proof
Since $n|N$ there exists $0\le\mu\le N-1$ so that $q=\go ^\mu.$ Using 
Proposition \ref{generalh} (2), define $A,X\in H^*$ as follows:
$$<A,a^ix^j>=\gd _{j,0}\go ^i\;\;\;<X,a^ix^j>=\gd _{j,1}$$
for all $0\le i\le N-1,0\le j\le r-1.$ Then it is straightforward to 
check that $A\in G(H^*)$ is of order $N,$ $X\in 
P_{A^\mu,\varepsilon}(H^*),$ $X^r=0,$ $XA=\go ^\nu AX,$ and finally that  
$H^*=k<A,X>=H_{\frac{N}{(N,\nu)},\go ^\nu,N,\mu}.$\qed 
\begin{Corollary}\label{new}
Let ${\cal H}_{\ga}={\cal H}_{n,q,N,\nu,\ga}\ne k\Z_N,$ and let $\go 
\in k$ be a primitive $Nth$ root of unity.
Then as an algebra ${\cal
H}_{\ga}^*\cong H_{\frac{n}{(n,\nu)},\go ^\nu,n,\mu}\#_\tau
k\Z_{N/n}$
is a crossed product of
a Radford Hopf algebra and a group algebra for some
$0\le\mu\le N-1$ so that $q=\go ^\mu.$ If moreover $\ga
\ne 0$ then
${\cal H}_{\ga}^*\cong H_{n,q,n,\nu}\#_\tau k\Z_{N/n}.$
\end{Corollary}
\proof
By the proof of Proposition \ref{generalh} (5), we have the following
short exact sequence of Hopf algebra maps:
$$k\Z_{N/n}\stackrel{i}{\hookrightarrow}{\cal H}_\ga 
\stackrel{\pi}{\raro}H_{n,q,n,\nu}.$$  
Dualizing yields the following short exact sequence of Hopf algebra
maps:
$$H_{n,q,n,\nu}^*\stackrel{\pi^*}{\hookrightarrow}{\cal H}_\ga ^*
\stackrel{i^*}{\raro}(k\Z_{N/n})^*.$$  
Since $(k\Z_{N/n})^*\cong k\Z_{N/n}$ as Hopf algebras 
and by Proposition \ref{rad}, $H_{n,q,n,\nu}^*\cong
H_{\frac{n}{(n,\nu)},\go ^\nu,n,\mu}$ for some
$0\le\mu\le N-1$ so that $q=\go ^\mu,$
it follows that
${\cal H}_{\ga}^*\cong H_{\frac{n}{(n,\nu)},\go ^\nu,n,\mu}\#_\tau
k\Z_{N/n}.$
If moreover $\ga \ne 0$ then $(n,\nu)=1$ and hence $\go =q
^{\nu^{-1}(mod\,n)}$ is a primitive $nth$ root of unity so that $\go ^\nu
=q.$ Thus, by Theorem \ref{selfdual}, $H_{n,q,n,\nu}^*\cong H_{n,q,n,\nu}$
and the result follows.\qed

We shall need the following lemma in which we find all the Hopf 
subalgebras of $H_{n,q,N,\nu}.$
\begin{Lemma}\label{subhopf}
Let $H=H_{n,q,N,\nu}\ne k\Z_N$ and $B\subseteq H$ be a Hopf 
subalgebra of $H.$ Then 
either $B=kG(B)$ where $G(B)$ is a subgroup of $(a),$ or there exists an 
integer $1\le s\le N$ so that $s|N,$ $(N/s)|\nu$ and
$B=k<a^{N/s},x>=H_{\frac{n}{(n,N/s)},q^{N/s}, 
s,\frac{\nu}{N/s}}.$ 
\end{Lemma}
\proof Suppose $B\ne kG(B).$ Then by Taft-Wilson Theorem [TW], $B$ must
contain a non-trivial skew primitive element $P.$ Without loss of 
generality we may assume that $P\in P_{a^i,1}(B)$ for some $0\le i\le 
N-1.$ By Proposition \ref{generalh}(3), $i=\nu,$ and hence $x\in B.$
 
Let $i:B\raro H$ be the inclusion map. Then $i^*:H^*\raro B^*$ is a 
surjection of Hopf algebras. By Proposition \ref{rad}, 
$H^*=k<A,X>=H_{\frac{N}{(N,\nu)},\go ^\nu, N,\mu}$ for some $0\le\mu\le 
N-1$ so that $q=\go ^\mu$ (where $A,X$ are as in the proof of Proposition 
\ref{rad}). Therefore, $B^*$ is generated as an algebra by 
$i^*(A),i^*(X)$ and $i^*(A)\in G(B^*),$ $i^*(X)\in P_{i^*(A)^{\mu},1}(B^*)$ 
(note that since $x\in B,$ $i^*(X)$ must be non-trivial). Set 
$s=|i^*(A)|.$ Then $s|N.$ It is now straightforward to check that 
there exists a surjection of Hopf algebras 
$H_{\frac{N}{(N,\nu)},\go ^\nu, s,\mu}\raro B^*=k<i^*(A),i^*(X)>$ (here 
$\mu$ means $\mu(mod\,s)$) which is also an injection by Proposition 
\ref{generalh}(3) and [M, Theorem 5.3.1]. Now, $\go ^{N/s}$ 
is a primitive $sth$ root of unity, and hence by Proposition \ref{rad}, 
$B=H_{\frac{N}{(N,\nu)},\go ^\nu, s,\mu}^*=
H_{\frac{s}{(s,\mu)},q^{N/s},s,\gb}$ for some $\gb$ so that $\go 
^{\nu}=\go ^{(N\gb )/s}.$ In particular the last equation implies that 
$(N/s)|\nu$ and $\gb =\frac{\nu }{N/s}(mod\,s).$ This concludes the 
proof of the lemma.\qed

It is known that $H_{(N,\nu,\go )}$ is quasitriangular if and
only if $N=2\nu$ and $\nu$ is odd [G1,R1]. In the following we
investigate the quasitriangularity of ${\cal H}_\ga .$
\begin{Theorem}\label{qt}
Let $k$ be an algebraically closed field of characteristic $0.$ Then 
${\cal H}_\ga ={\cal H}_{n,q,N,\nu,\ga }$ is quasitriangular if 
and only if either ${\cal H}_{\ga}=k\Z_N$ or ${\cal H}_{\ga} 
=H_{(2\nu,\nu,\go )},$ where $\nu$ is odd. 
\end{Theorem}
\proof
We have to prove the 'only if' part. We first show that if ${\cal H}_\ga $ is
quasitriangular then $\ga =0.$
Suppose on the contrary that $({\cal H}_\ga ,R)$ is quasitriangular and
$\ga \ne 0.$ Write
$R=\sum R^{(1)}\ot R^{(2)}$ in the shortest possible way. Then
$A=sp\{R^{(1)}\}$ and $B=sp\{R^{(2)}\}$ are Hopf subalgebras of ${\cal
H}_\ga$ and $A\cong B^{*cop}.$ In particular $A,$ $B,$ $A^*$ and $B^*$ are
pointed. Suppose $A\ne kG(A).$ Then using similar arguments to those used in 
the first paragraph of the proof of Lemma \ref{subhopf}
yields that $a^{\nu },x\in A.$
Therefore, we have an inclusion of Hopf algebras $k<a^{\nu },x>\subseteq 
A,$ and hence a surjection of Hopf algebras $A^*\raro k<a^{\nu },x>^*.$ 
Since $A^*$ is pointed so is $k<a^{\nu },x>^*.$ But, $k<a^{\nu 
},x>^*={\cal H}_{n,q^{\nu },\frac{N}{(N,\nu)},1,\ga }$ is {\em not} 
pointed by Theorem \ref{dualpoint}, which is a contradiction. We conclude that
$A=kG(A)=k(a^m)$ for some $m|N,$ and since $B\cong A^{*cop}$ it follows
that $B=k(a^m)$ as well. Therefore, $R\in k(a)\ot k(a).$ Let $\go \in k$
be a primitive $Nth$ root of unity. By [R1, page 219], there exists $0\le
l\le N-1$ so that 
\begin{equation}\label{r}
R=\frac{1}{N}\left(\sum_{i,j=0}^{N-1}\go ^{-ij}a^i\ot a^{jl}\right).
\end{equation} 
By (QT.5), $(a^\nu \ot x+x\ot 1)R=R(x\ot a^\nu +1\ot x)$; that is, 
\begin{equation}\label{eq}
\frac{1}{N}\left(\sum_{i,j=0}^{N-1}\go ^{-ij}q^{jl}a^{i+\nu}\ot
a^{jl}x\right)+\frac{1}{N}\left(\sum_{i,j=0}^{N-1}\go
^{-ij}q^ia^ix\ot a^{jl}\right)=
\end{equation}
$$\frac{1}{N}\left(\sum_{i,j=0}^{N-1}\go ^{-ij}a^ix\ot
a^{jl+\nu}\right)+\frac{1}{N}\left(\sum_{i,j=0}^{N-1}\go
^{-ij}a^i\ot a^{jl}x\right).$$
Thus, $\frac{1}{N}\go ^{-ij}q^{jl}=\frac{1}{N}\go ^{-(i+\nu )j},$
i.e. $q^{jl}=\go ^{-\nu j}$ for all $j.$ In particular, for $j=n$ we
obtain $1=q^{nl}=\go ^{-\nu n}.$ Thus, $N|\nu n$ which is a contradiction.
Therefore $\ga =0.$

Suppose $H_{n,q,N,\nu}\ne k\Z_N,$ $(H_{n,q,N,\nu},R)$ is 
quasitriangular, and let $A$ and $B$ 
be defined as before. Suppose $A\ne kG(A)$ (equivalently, $B\ne kG(B)$  
since $A\cong B^{*cop}$). By Lemma \ref{subhopf}, there exist integers 
$1\le s,s'\le N$ so that $s,s'|N,$ $(N/s),(N/s')|\nu$ and
$$A=k<a^{N/s},x>=H_{\frac{n}{(n,N/s)},q^{N/s},s,
\frac{\nu}{N/s}}$$
and
$$B=k<a^{N/s'},y>=H_{\frac{n}{(n,N/s')},q^{N/s'},s',  
\frac{\nu}{N/s'}}.$$
Since $A\cong B^{*cop},$ it follows from Proposition \ref{rad} that 
$s=s',$ and hence by [G2, Lemma 1.1.2], there exists a primitive
$sth$ root of unity $\chi $ so that $q^{N/s}=\chi
^{\frac{\nu}{N/s}}$ and $q^{N/s}=\chi
^{-\frac{\nu}{N/s}}.$ Therefore $N|2\nu,$ and since $\nu <N,$ we conclude 
that $N=2\nu.$ Now, since $A$ and $B$ generate $H_{n,q,N,\nu}$
as a Hopf algebra, we must have $s=N.$ Therefore we may assume that $\chi=\go$ 
and hence
$q=\go ^{\nu}.$ Since by assumption $H_{n,q,N,\nu}\ne k\Z_N,$ we have that 
$1<r=|q^{\nu}|=|\go^{\nu^2}|,$ and hence $2\nu\Nb \nu^2$ which 
implies that $\nu$ is odd. Therefore, 
$A=B=H_{n,q,N,\nu}=H_{(2\nu,\nu,\go)}$ where $\nu$ is odd.

If $A=B=k(a^m)$ then again $R$ must be of the form (\ref{r}). But then it
follows from (\ref{eq}) that $l$ must be invertible in $\Z_N,$ and hence
that on one hand $q=\go ^{-\nu l^{-1}}$ and on the other hand, $q=\go
^{\nu l^{-1}}.$ Therefore, $N=2\nu$ and $q=\go ^{\nu}=-1.$ By [R1, 
Corollary 3(a)], $H_{n,q,N,\nu}=H_{(2\nu,\nu,\go)}$ where $\nu$ is odd.   
This completes the proof of the theorem.\qed  

We conclude this section with an example.
\begin{Example}\label{8}
{\rm The least dimensional new Hopf algebra in our family is an 
$8-$dimensional
Hopf algebra, namely $H_8={\cal H}_{2,-1,4,1,1},$ generated as an algebra
by $a\in G(H_8)$ and $x\in P_{a,1}(H_8)$ which satisfy the relations 
$$a^4=1,\;x^2=a^2-1\;\;and\;\;xa=-ax.$$
By Proposition \ref{generalh}, as an algebra $H_8=k\Z_2\#_\sigma H_4$ 
is an extension of $k\Z_2$ by Sweedler's Hopf algebra $H_4.$ 
By Theorem \ref{qt}, $H_8$ is not quasitriangular even though $k\Z_2$ and 
$H_4$ are. By Lemma \ref{dualgl},
$G(H_8^*)=\{1,\phi\}$ is of order $2,$ where $\phi:H_8\raro k$ is
determined by $\phi(a)=-1$ and $\phi(x)=0.$ 
By Corollary \ref{new}, as an algebra $H_8^*=H_4\#_\tau k\Z_2$ is an 
extension of $H_4$ by $k\Z_2.$ By Theorem \ref{dualpoint}, $H_8^*$ is
not pointed and it contains a (unique) $4-$dimensional simple sub
coalgebra. In particular $H_8$ is not self dual.} 
\end{Example}
\section{The Family ${\cal U}_{(n,N,\nu,q,\ga ,\gb ,\gga )}$}
In this section we first construct a new family of finite dimensional 
pointed and unimodular Hopf algebras, denoted by ${\cal U}_{(n,N,\nu,q
,\ga ,\gb ,\gga )},$ 
which generalizes Radford's family $U_{(N,\nu,\go)}.$ Second, we 
show that over any infinite field which contains a primitive $nth$ root
of unity, our new family contains infinitely many non-isomorphic Hopf
algebras of any dimension of the form $Nn^2,$ where $2<n<N$ are 
integers so that $n$ divides $N.$ Thus, we prove that Kaplansky's 
10th conjecture [K] is false.
Third, we use
Theorem \ref{dualpoint} to 
generalize [G3, Theorem 1.2.2] and show that over an algebraically
closed field of characteristic zero the family $U_{(N,\nu,\go)},$ any $N,$
is unique with respect to properties (U.1)-(U.3).

We start by introducing ${\cal U}_{(n,N,\nu,q,\ga ,\gb ,\gga )}.$
Let $n>1,$ $\nu$ and $N$ be positive integers such that $1\leq\nu <N$ and
$n|N,$ and let $\ga ,\gb ,\gga \in k.$ Suppose that $q\in k$ is a 
primitive $nth$ root of unity, and let $r=|q^\nu|.$ For $\ga
=\gb =0$ we define ${\cal U}_{(n,N,\nu,q,0,0,\gga )}$ exactly as
$U_{(N,\nu,\go)}$ except
that now $yx-q^{-\nu}xy=\gga (a^{2\nu}-1).$ 
For $(\ga ,\gb )\ne (0,0)$ we assume that $(n,\nu)=1$ 
(hence $r=n$), and $N\Nb\nu n.$ In this case, as an algebra 
${\cal U}_{\ga ,\gb ,\gga }={\cal U}_{(n,N,\nu,q,\ga ,\gb ,\gga )}$ is 
generated by $a,x$ and $y$ which satisfy the relations
$$a^N=1,\;x^n=\ga (a^{\nu n}-1),\;y^n=\gb (a^{\nu 
n}-1),$$ 
$$xa=qax,\;ya=q\minus 
ay\;\;and\;\;yx-q^{-\nu}xy=\gga (a^{2\nu}-1).$$
The coalgebra structure of ${\cal U}_{\ga ,\gb ,\gga }$ is determined by
$$\gD(a)=a\ot a,\;\;\gD(x)=x\ot a^\nu+1\ot x,\;\;\gD(y)=y\ot a^\nu+1\ot y,$$
$$\gep(a)=1,\;\;\;\gep(x)=0 \;\;\;{\rm and}\;\;\;\gep(y)=0.$$
The antipode of ${\cal U}_{\ga ,\gb ,\gga }$ is determined by
$$s(a)=a\minus,\;\;\;s(x)=-q^{-\nu}a^{-\nu}x\;\;\;{\rm and}\;\;\;
s(y)=-q^\nu a^{-\nu} y.$$
Note that if $\go$ is a primitive $Nth$ root of unity and $N\Nb\nu^2$, then
$U_{(N,\nu,\go)}={\cal U}_{(\frac{N}{(N,\nu )},N,\nu,\go ^\nu,0,0,1)},$ and
${\cal U}_{(\frac{N}{(N,\nu )},N,\nu,\go ^\nu,0,0,\gga )}\cong 
U_{(N,\nu,\go)}$ as 
Hopf algebras for any $\gga \in k^*$ which has a square root in $k.$
\begin{Remark}\label{radf2}
{\rm Suppose $(\ga ,\gb )\ne (0,0).$ Then to construct ${\cal 
U}_{(n,N,\nu,q,\ga ,\gb ,\gga )}$ we follow
the line of argument and notations used by Radford in order to construct
$U_{(N,\nu,\go)}$ [R1, 5.2] (see also Remark \ref{radf1}). We let $C$ be the 
$5-$dimensional coalgebra
over $k$ with basis $\{A,B,D,X,Y\}$ whose structure is determined by
$$\gD(A)=A\ot A,\;\;\gD(B)=B\ot B,\;\;\gD(D)=D\ot D,$$
$$\gD(X)=X\ot B+ D\ot X\;\;and\;\;\gD(Y)=Y\ot B+ D\ot Y,$$
and let $I$ be the ideal of the tensor algebra $T(C)$ generated by
$$A^N-1,\;\;B-A^{\nu},\;\;D-1,\;\;XA-qAX,\;\;YA-q^{-1}AY,$$
$$X^n-\ga (A^{\nu n}-1),\;\;Y^n-\gb (A^{\nu 
n}-1)\;\;and\;\;YX-q^{-\nu}XY-\gga (A^{2\nu}-1).$$
The details are straightforward to work out.}
\end{Remark}
\begin{Proposition}\label{generalu}
Let ${\cal U}_{\ga ,\gb ,\gga }={\cal U}_{(n,N,\nu,q,\ga ,\gb ,\gga  )}.$ 
Then: 
\ben
\item ${\cal U}_{\ga ,\gb ,\gga }$ is a pointed Hopf algebra with $G({\cal
      U}_{\ga ,\gb ,\gga })=(a)$ a cyclic group of order $N.$
\item The set $\{a^ix^jy^l|0\le i \le N-1,0\le j,l \le r-1\}$ forms a linear 
      basis of ${\cal U}_{\ga ,\gb ,\gga },$ hence $dim{\cal U}_{\ga ,\gb
      ,\gga }=Nr^2.$
\item $P_{a^\nu,1}({\cal U}_{\ga ,\gb ,\gga })=sp_k\{a^\nu-1,x,y\},$
      $P_{a^v,1}({\cal U}_{\ga ,\gb ,\gga })=sp_k\bbraces{a^v-1}$ for
      all $v\neq\nu,$ and $({\cal U}_{\ga ,\gb ,\gga 
      })_1=sp_k\{a^i,a^jx,a^ly|0\le i,j,l\le N-1\}.$
\item The element $\lambda=\left (\frac{1}{N}\sum_{i=0}^{N-1}a^i\right 
      )x^{r-1}y^{r-1}$ is a non-zero 2-sided integral of ${\cal U}_{\ga 
      ,\gb ,\gga }.$ In particular ${\cal U}_{\ga ,\gb ,\gga }$ is 
      unimodular. 
\item If $\ga ,\gb ,\gga \ne 0$ and $\sqrt\gga \in k,$ then as an
      algebra, ${\cal U}_{\ga ,\gb
      ,\gga }\cong k\Z_{N/n}\#_{\sigma}U_{(n,\nu,\go)}$ is a crossed
      product of a group algebra and a Radford Hopf algebra for some 
      primitive $nth$ root of unity $\go .$ 
\item Let $\ga, \gb ,\gga,\ga ',\gb ',\gga 
'\in k^*.$ Then: \ben
\item ${\cal U}_{\ga ,\gb ,0}\cong {\cal U}_{\ga ',\gb ',0}$ as Hopf 
      algebras.
\item Suppose $N\Nb 2\nu.$ Then ${\cal U}_{\ga ,\gb ,\gga }\cong {\cal 
      U}_{\ga ',\gb ',\gga '}$ as Hopf algebras if and only if 
      $\frac{\ga \gb }{\gga ^n}=\frac{\ga '\gb '}{\gga '^n}.$
\een
\een
\end{Proposition}
\proof
1. Follows from Lemma \ref{point}.\\
2. For $U_{(N,\nu,\go)}$ it was proved in [R1, Proposition 10 b)] using
Bergman's Diamond Lemma. The proof for ${\cal U}_{\ga ,\gb ,\gga }$ is 
very similar, and we merely 
sketch it following the line of argument from the proof of [R1,
Proposition 10 b)]. Using the notation of Remark \ref{radf2},
let $T(C)=k\{A,B,D,X,Y\}$ where $A<B<D<X<Y$ and the substitution rules are:
$$A^N\leftarrow 1,\;\;B\leftarrow A^{\nu},\;\;D\leftarrow 1,\;\;XA\leftarrow
qAX,\;\;YA\leftarrow q^{-1}AY,$$
$$X^n \leftarrow \ga (A^{\nu n}-1),\;\;Y^n \leftarrow \gb (A^{\nu 
n}-1)\;\;and\;\;YX\leftarrow q^{-\nu}XY+\gga (A^{2\nu}-1).$$
The only ambiguities to consider are overlap ambiguities, and they are all
easily resolved. We just note that in the case $(\ga ,\gb )\ne (0,0)$ it is
necessary to have $r=n$ (so that $q^r=q^n=1$) when resolving, for instance, 
the overlap ambiguity $(YX)X^{n-1}=Y(XX^{n-1}).$\\ 
3. The proof follows the line of argument and notations of the proof of 
[G2, Lemma 1.2.1]. Specifically, 
let $A=k(a),$ $B=k<a,x>,$ $ H=k<a,y>,$ $\{B_m\}$ be the coradical filtration
of $B$ and $\{H_m\}$ be the coradical filtration of $H.$ Note that
$B_0=H_0=A.$ Set $$ B(1)=Ax,\;\;\;\;H(1)=Ay,$$
$$ B(m)=B(1)^m,\;\;\;\;H(m)=H(1)^m\;\;\;for\;\;\;1\leq m$$
$$and\;\;\;B(0)=H(0)=A.$$
Using parts 2 and 3 of Proposition 2.1 it is straightforward to verify that:
$$B=\Op_{i=0}^{r-1}B(i),\;\;\;B(0)=B_0,\;\;\;B_1=B(0)\op B(1),$$
$$H=\Op_{i=0}^{r-1}H(i),\;\;\;H(0)=H_0\;\;\;and\;\;\;H_1=H(0)\op H(1).$$
Therefore, $B$ and $H$ are coradically graded [CM, Section 2].
Since ${\cal U}_{\ga ,\gb ,\gga }=B\ot_{A}H$ it follows from [CM, Lemma 
2.3] that ${\cal U}_{\ga ,\gb ,\gga }$ is coradically graded where,
$$ {\cal U}_{\ga ,\gb ,\gga }(m)=\sum_{i=0}^{m}B(i)\ot_{A}H(m-i).$$
In particular, $({\cal U}_{\ga ,\gb ,\gga })_0={\cal U}_{\ga ,\gb 
,\gga }(0)=A$ and $${\cal 
U}_{\ga ,\gb ,\gga }(1)=B(0)\ot_{A}H(1)+B(1)\ot_{A}H(0)=H(1)+B(1).$$
Thus $({\cal U}_{\ga ,\gb ,\gga })_1={\cal U}_{\ga ,\gb ,\gga }(0)\op 
{\cal U}_{\ga ,\gb ,\gga }(1)=A\op Ax\op Ay$ and the result follows.\\
4. Since $\frac{1}{N}\left(\sum_{i=0}^{N-1}a^i\right)$ is a 2-sided
integral of $k(a),$ and $a$ commutes with $x^{r-1}y^{r-1}$ it follows
that $a\lambda=\lambda a=\lambda=\varepsilon(a)\lambda.$ By
the same reason, since at any rate $x^r,y^r\in k(a)\cap ker(\varepsilon)$
it follows that
$x\lambda=\lambda y=0$ (see the proof of Proposition \ref{generalh}
(4)). We show now that $\lambda x=0.$ Indeed,
by simple induction we have 
$$y^ix=q^{-i\nu}xy^i+
\gga\left(q^{-(i-1)\nu}\sum_{j=0}^{i-1}q^{-j\nu}\right)a^{2\nu}y^{i-1}-
\gga\left(\sum_{j=0}^{i-1}q^{-j\nu}\right)y^{i-1}.$$ 
Using this, and the fact that $x^r\in k(a)\cap ker(\varepsilon)$ we compute: 
\begin{eqnarray*}
\lefteqn{\lambda
x=\left(\left(\frac{1}{N}\sum_{i=0}^{N-1}a^i\right)x^{r-1}y^{r-1}\right)x}\\
& &=\left(\frac{1}{N}\sum_{i=0}^{N-1}a^i\right)x^{r-1}\left(q^\nu
xy^{r-1}-\gga q^{3\nu}a^{2\nu}y^{r-2}+\gga q^\nu y^{r-2}\right)\\
& &=\gga \left(\frac{1}{N}\sum_{i=0}^{N-1}a^i\right)x^{r-1}
\left(-q^{3\nu}a^{2\nu}+q^\nu 1\right)y^{r-2}\\
& &=\gga \left(\frac{1}{N}\sum_{i=0}^{N-1}a^i\right)\left(-q^{\nu}
a^{2\nu}+q^\nu 1\right)x^{r-1}y^{r-2}\\
& &=\varepsilon\left(-q^{\nu}a^{2\nu}+q^\nu
1\right)\gga \left(\frac{1}{N}\sum_{i=0}^{N-1}a^i\right)x^{r-1}y^{r-2}=
0=\varepsilon(x)\lambda.
\end{eqnarray*}
Similarly, $y\lambda =0,$ and the result follows.\\
5. Since $a^n$ is central in ${\cal U}_{\ga ,\gb ,\gga }$ it follows that
$K=k(a^n)\cong k\Z_{N/n}$
is a central Hopf subalgebra of ${\cal U}_{\ga ,\gb ,\gga },$ hence is
normal. Since $(n,\nu)=1$ there exist integers $a,b$ such that 
$1=an+b\nu.$ Set $\go =q^b.$ Then $\go $ is a primitive $nth$ root of 
unity, and $\go ^{\nu}=q.$ Let $\pi :{\cal U}_{\ga ,\gb ,\gga } \raro 
U_{(n,\nu,\go)}$ be the Hopf algebra
surjection determined by $\pi (a)=a,$ $\pi(x)=\sqrt\gga x$ and $\pi(y)= 
\sqrt\gga y.$ Then ${\cal U}_{\ga ,\gb ,\gga }K^+=ker\pi,$ hence
the following is a short exact sequence of Hopf algebra maps:
$$K\stackrel{i}{\hookrightarrow}{\cal U}_{\ga ,\gb ,\gga } 
\stackrel{\pi}{\raro}U_{(n,\nu,\go)}$$
where $i$ is the inclusion map, and the result follows.\\
6. First note that since the ground field contains a primitive $Nth$ 
root of unity and $n$ divides $N,$ it also contains a primitive 
$nth$ root of unity. Now, following the lines of the proof of Proposition
\ref{generalh} (6), it is straightforward to check, using part 3 of 
this proposition, that if $f:{\cal
U}_{\ga ,\gb ,\gga }=k<a,x,y>\raro {\cal U}_{\ga ',\gb ',\gga '}=k<g,z,w>$
is an isomorphism of Hopf algebras then $f(a)=g,$ $f(x)=u(\frac{\ga}{\ga
'})^{1/n}z$ and $f(y)=v(\frac{\gb }{\gb '})^{1/n}w$ for some $nth$ roots 
of unity $u,v\in k.$ If $\gga =\gga '=0$ 
then such an $f$ is indeed an isomorphism of Hopf algebras. If $\gga ,
\gga '\ne 0$ then since $N\Nb 2\nu$
it follows that $g^{2\nu}-1\ne 0,$ and hence that
$f(yx-q^{-\nu}xy)=\gga f(a^{2\nu}-1)$ if and only if $uv(\frac{\ga \gb }{\ga
'\gb '})^{1/n}(zw-q^{-\nu}wz)=\gga (g^{2\nu}-1)$ if and only if
$u^{-1}v^{-1}(\frac{\ga \gb }{\ga '\gb '})^{-1/n}\gga =\gga 
'$ if and only if $\frac{\ga \gb }{\gga ^n}=\frac{\ga '\gb '}{\gga '^n}.$
\qed
\begin{Corollary}\label{kap}
Kaplansky's 10th conjecture is false over any infinite 
field which contains a primitive $nth$ root of unity for some $n>2.$ 
\end{Corollary}
\proof
The conjecture was that there exist only finitely many non-isomorphic Hopf 
algebras of any given dimension. But, since the field is infinite, it 
follows by Proposition \ref{generalu} (6) that the family ${\cal
U}_{\ga ,\gb ,\gga },$ where $\ga ,\gb ,\gga \ne 0,$ contains infinitely 
many non-isomorphic Hopf algebras of any dimension of the form  
$Nn^2,$ where $2<n<N$ are integers so that $n$ divides $N.\qed$

In what follows we characterize Radford's $U_{(N,\nu,\go)}$ as a 
distinguished sub-family of our new family via minimal quasitriangularity. 
\begin{Lemma}\label{1}
Let $k$ be an  algebraically closed field of characteristic zero.
Let $U$ be a Hopf algebra over $k$ which satisfies the following:
\ben
\item $G(U)=(a)$ is a cyclic group  of order $N.$
\item There exists an integer $1\leq\nu< N$ so that 
      $\dim_k P_{a^\nu , 1}(U)=3$.
\een
Then $U$ contains two Hopf subalgebras, ${\cal H}_{n_1, \go_1,
N, \nu ,\ga _1}=k<a,x>$ and ${\cal H}_{n_2,\go_2,N,\nu ,\ga _2}=k<a,y>$
for some $n_1,n_2\in\NN$, $\ga _1,\ga _2,\go_1,\go_2\in
k.$ Moreover, ${\cal H}_{n_1, \go_1, N, \nu ,\ga _1}\cap
{\cal H}_{n_2,\go_2,N,\nu ,\ga _2}=kG(U).$
\end{Lemma}
\proof Since  $kG(U)$ is a commutative semisimple  algebra
over an algebraically  closed field of characteristic zero 
all $kG(U)-$modules are direct sums of $1-$dimensional $kG(U)-$
submodules. Since $kG(U)$ acts on $P_{a^\nu,1}(U)$ by conjugation it
follows that $P_{a^\nu,1}(U)$ is a direct sum of 
three $1-$dimensional $kG(U)-$submodules. Therefore we may assume that 
there exists a basis $\{a^\nu-1,x,y\}$ for $P_{a^\nu,1}(U)$ so that    
$$ x a=\go_1 ax \;\; and \;\; ya=\go_2 ay$$
where $\go_i$ is an $Nth$ root of unity  for $i=1,2$.
Set $n_i=|\go_i|$ and $r_i=|\go_i^\nu|$ for $i=1,2$.
Following the line of argument of the proof of Theorem \ref{h}
yields that $x^{r_1}=\ga _1(a^{\nu r_1} -1)$ and $y^{r_2}=\ga _2(a^{\nu 
r_2} -1)$ for some $\ga _1,\ga _2\in k.$ Therefore there exist surjections 
of Hopf algebras ${\cal H}_{n_1,\go_1,N,\nu ,\ga _1}\raro k<a,x>$ and 
${\cal H}_{n_2,\go_2,N,\nu ,\ga _2}\raro k<a,y>$ which are also 
injections by Proposition \ref{generalh}(3) and [M, Theorem 5.3.1]. 
This 
completes the proof of the lemma.\qed \begin{Proposition}\label{odd}
Let $U$ be a finite dimensional Hopf algebra over an algebraically
closed field $k$ of characteristic zero. Then the following are equivalent:
\ben
\item $U$ satisfies the following conditions:
\begin{description}
\item[(i)] There exist integers $1\leq \nu< N,$ $N\Nb \nu^2,$  
           so that as an algebra $U$ is generated by a grouplike element 
           $a,$ 
           of order $N,$ and by the space $P_{a^\nu,1}(U)$ of dimension $3.$
\item[(ii)] $U\cong U^{*\;cop}$ as Hopf algebras.
\end{description}
\item $U=U_{(2\nu,\nu,\go)}$, where $\nu$ is odd.
\een
\end{Proposition}
\proof  Suppose 2 holds. Then (i) is just (U.1). We prove (ii).
By Proposition \ref{generalu}(2), the set $\{a^i x^j y^l|0\leq
i< N,0\leq j,l<r \}$ forms a linear basis of $U$.
Define  $A,X$ and $Y$ in $U^*$ as follows:
$$<A,a^i x^jy^l>=\go^i\gd_{j,0}\gd_{l,0},\;\;\;<X,a^i x^jy^l>=\go^{\nu i}
\gd_{j,1}\gd_{l,0}= (-1)^i \gd_{j,1}\gd_{l,0}$$ 
and
$$<Y,a^ix^jy^l>=\go^{\nu i}\gd_{j,0} \gd_{l,1}=(-1)^i\gd_{j,0}\gd_{l,1}$$
for $\;0\leq i\leq 2\nu -1$ and $0\leq j,l\leq r-1$.
Using a similar argument to the one used in the proof of [G2, Lemma 1.1.2]
yields that the set $\{A^i X^j Y^l|0\leq i< N,0\leq j,l<r \}$ forms a
linear basis  of $U^*$. It is now straightforward to check that the  map 
$f:U\raro U^{*\;cop}$ given by
$$f(a)=A,\;\; f(x)=X\;\;{\rm and}\;\; f(y)=Y$$
is an isomorphism of Hopf algebras.

Assume 1 holds. It follows  from Lemma \ref{1} that there exist 
$x,y\in\Pa$ such that $k<a,x>={\cal H}_{n_1,\grr,N,\nu ,\ga _1}$  and
$k<a,y>={\cal H}_{n_2,\go_2,N,\nu ,\ga _2}.$ Let $f:U\raro U^{*\;cop}$ be
the isomorphism of Hopf algebras which
exists by (ii).
Set $f(a)=A,\;\;f(x)=X$ and $f(y)=Y$. It follows that $A\in G(U^*)$
and  $\go =<A,a>\in k$ is a primitive $Nth$ root of unity.
By Lemma \ref{1} there exists an inclusion of Hopf algebras, $i:
{\cal H}_{n_1,\grr,N,\nu ,\ga _1}\raro U$, hence there exists a surjection
of Hopf algebras,
$i^*:U^*\raro ({\cal H}_{n_1,\grr,N,\nu ,\ga _1})^*$. Since $i^*(A),\;\;
i^*(X)$ and
$i^*(Y)$ generate  $({\cal H}_{n_1,\grr,N,\nu ,\ga _1})^*$ as an algebra
it follows
that either $i^*(X)\ne 0$ or $i^*(Y)\ne 0$. Therefore using a similar 
argument to the one 
used in the third paragraph of the proof of [G2, Lemma 1.1.2] we
conclude that either $<X,x>\neq 0$ or 
$<Y,x>\neq 0$. Let $Z$ be $X$ or $Y$ so that $<Z,x>\neq 0$. 
Then $\grr<Z,ax>=<Z,xa>$ yields
$\grr=\go^{-\nu}$. In the same manner either  $<X,y>\neq 0$ or $<Y,y>\neq
0$ and the equality $ya=\go_2 ay$ yields $\go_2 =\go^{-\nu}.$ In
particular $\go_1=\go_2$. Therefore $<ZA,x>=<\go^{-
\nu}AZ,x>$ implies $\go^{-\nu}=\go^\nu$. Thus, $\go^{2\nu}=1$ and
$N|2\nu$. Since $\nu<N$ we have $N=2\nu$. But $N\Nb\nu^2$ hence we
conclude that $\nu$ is odd. Moreover, since $\go^\nu=-1$ it follows from the
definition of $\gD(x)$ and $\gD(y)$ that $\gD(yx+xy)=(xy+yx)\ot 1+1\ot 
(xy+yx)$. By Lemma \ref{trivial}, $xy+yx=0.$ Finally, $r=|\go ^{\nu 
^2}|=\frac{2\nu }{(2\nu ,\nu ^2)}=2,$
hence $x^2=\ga _1(a^{2\nu }-1)=0.$ Similarly, $y^2=0.$ We conclude that 
there exists a surjection of Hopf algebras $U_{(2\nu,\nu,\go)}\raro U,$ 
which is also an injection by (U.2) and [M, Theorem 5.3.1]. This completes 
the proof of the proposition.\qed 

We are ready now to characterize Radford's $U_{(N,\nu,\go)}$ via minimal
quasitriangularity.
\begin{Theorem}\label{u}
Let $k$ be an algebraically closed field of characteristic zero. Let $U$
be a Hopf algebra which satisfies the following:
\ben
\item There exist integers $1\leq \nu< N,$ $N\Nb \nu^2,$
      so that as an algebra $U$ is generated by a grouplike element $a,$
      of order $N,$ and by the space $P_{a^\nu,1}(U)$ of dimension $3.$
\item For any positive integer $i\neq \nu(mod\,N),$ $dim_k P_{
      a^i,1}(U)=1.$
\item $U$ is minimal quasitriangular.
\een
Then there exists a primitive $Nth$ root of unity $\go\in k$ such that 
$U=U_{(N,\nu,\go)}$.
\end{Theorem}
\proof
By Lemma \ref{1}, there exist $x,y\in\Pa$ so that the  
set $\{a^\nu -1,x,y\}$ forms a linear basis for $\Pa$ and the following
relations are satisfied: 
$$a^N=1,\;xa=\go_2 ax,\;ya=\grr ay,\;x^{r_2}=\ga _2(a^{\nu 
r_2}-1)\;and\;y^{r_1}=\ga _1(a^{\nu r_1}-1)$$
where $\go_i$ is a primitive $n_ith$ root of unity  
and $r_i=|\go_i^\nu|$ for $i=1,2.$
Moreover, as an algebra $U$ is generated by $a,x$ and $y$. 
Since $U$ is minimal quasitriangular it follows by [R2, Proposition 2(d)
and Theorem 1]
that there exist two Hopf subalgebras $B$ and $H$ of $U$ such that $U$ is 
generated by $B$ and $H$ and $B^{*\;cop}\cong H$ as Hopf algebras. In
particular  $B\neq kG(B)$ and $H\neq kG(H),$ otherwise
$B=H=U\subseteq k(a)$ which is a contradiction. Note also that $B^*$ and
$H^*$ are necessarily pointed Hopf algebras.

Let $C\cntd U$ be a Hopf
subalgebra such that $C\neq kG(C)$ and $C^*$ is pointed. Since $U$ is pointed
$C$ is pointed as
well. Since $C\neq kG(C)$ it follows from Taft-Wilson Theorem [TW] that
$C$ contains
a non-trivial skew primitive element $P$. Without loss of generality we
may assume that $P\in P_{a^i,1}(U)$ for some  $0\leq i\leq N-1$.
By assumption 2, $i=\nu$. Thus $\gD(P)=P\ot a^\nu +1\ot P$
and $k(a^\nu)\sub C$. Since $P$ must be  a linear  combination of $x,y$
and $a^\nu -1$ we conclude, using Lemma \ref{subhopf}, that there exists an 
integer 
$1\leq s \leq N$ such that $s|N,$ $(N/s)|\nu$ and $C$ is
one of the following distinct Hopf subalgebras of $U$:
$$k<a^{N/s},x,y>,\;\; k<a^{N/s},x>\;\;{\rm or}\;\;k<a^{N/s},y>.$$
If in addition $\grr^{N/s}=\go_2^{N/s}$ then $C$ may also equal
$$k<a^{N/s},\ga x+\gb y>$$
where $\ga,\gb$ are fixed elements in $k^*.$
Since $B$ and $H$ generate $U$ it follows by the above that the mutually exclusive
possibilities for $B$ and $H$ are the following:
\ben
 \item $B=H=U.$ 
 \item $B=k<a,x>$ and $H=k<a,y>$ (or vice versa).
 \item $B=k<a,x>$ and $H=k<a,\ga x+\gb y>$ (or vice 
       versa) where $\ga,\gb$ are fixed elements in $k^*.$
 \item $B=k<a,y>$ and $H=k<a,\ga x+\gb y>$ (or vice 
       versa) where $\ga,\gb$ are fixed elements in $k^*.$
 \item $B=k<a, \ga x+\gb y>$ and $H=k<a,\gga x+\gd y>$ where 
       $\ga,\gb,\gga ,\gd $ are fixed elements in $k^*,$
       so that $\ga\gd-\gb\gamma\neq 0$. 
\een

If 1 holds then by Proposition \ref{odd}, $U=U_{(2\nu,\nu,\go)}$
and we are done. 

If 3-5 hold then $\grr =\go _2.$ 
Since $B^*$ and $H^*$ are pointed it follows from Theorem
\ref{dualpoint} that for all 
$\ga,\gb,\gga,\gd \in k$ such that $\ga\gb\neq 0$  and $\gga\gd\neq 0$
$$k<a,\ga x+\gb y>\cong k<a,\gga x+\gd y>\cong H_{n_1,\grr,N,\nu}$$ 
as Hopf algebras. Therefore since $B^{*\;cop}\cong H$ as Hopf algebras it
follows from [G2, Lemma 1.1.2] that there exists a
primitive $Nth$ root of unity  $\go$, such that
$\go^{-\nu}=\grr=\go_2=\go^\nu$.
Therefore $N=2\nu,\;\nu$ is odd and $xy+yx=0$ (see the proof of Proposition
\ref{odd}). Thus $U=U_{(2\nu,\nu,\go)}$ and we are done.

Suppose 2 holds. Since $B^*$ and $H^*$ are pointed it follows from Theorem 
\ref{dualpoint} that $B=H_{n_1,\grr,N,\nu}$ and $H=H_{n_2,\go _2,N,\nu}.$
Thus, it follows from [G2, Lemma 1.1.2]
that there exists a primitive $Nth$ root of unity $\go$, such that
$\grr=\go^{-\nu}$ and $\go_2=\go^{\nu}$. In particular $r_1=r_2=r,$
so $x^r=y^r=0.$ By [R2, Theorem 2], there exists a surjection of Hopf 
algebras  $\pi :D(B)\raro U$ such that \mat{\pi_{\mid_B}= id_B} and 
\mat{\pi_{\mid_{B}*cop}:B^{*cop}{\raro}H} is an isomorphism of Hopf 
algebras. Let $A$ and $X$ be the canonical generators of $B^*.$ Then
using [G2, Lemma 1.1.2] and [R1, Theorem 4 d)] we conclude
that $\pi$ is given by: 
$$\pi(\gep\bab a)=a,\;\;\pi(\gep\bab x)=x,\;\;\pi(A\bab 1)=a^m\;\;{\rm  
and}\;\;  \pi(X\bab 1) =\gb y$$ 
for some $\gb\in k^*,$ and an integer $m$ 
which satisfies $(m,N)=1$ and $m\nu = \nu (mod\; N)$.
Since the relation
$$(\gep\bab x)(X\bab 1)=\go^{\nu^2}(A^\nu \bab a^\nu+X\bab x-\gep \bab 1)$$
holds in $D(B)$ we conclude that the relation $$-\gb yx+\gb\go^{-\nu^2} 
xy=a^{2\nu}-1$$ 
holds in $U$. Hence after replacing $(-\gb y)$ by $y$ we get that 
$U=U_{(N,\nu,\go )}$. This completes the proof of the theorem.\qed

\end{document}